\documentclass[12pt]{article}
\usepackage{amsmath}
\usepackage{amssymb}
\usepackage{vatola}

\def\q{\quad}
\def\qq{\qquad}

\def\mod#1{\ (\text{\rm mod}\ #1)}
\def\t{\text}
\def\f{\frac}
\def\e{\equiv}
\def\b{\binom}
\def\sls#1#2{(\f{#1}{#2})}
 \def\ls#1#2{\big(\f{#1}{#2}\big)}
\def\Ls#1#2{\Big(\f{#1}{#2}\Big)}

\let \pro=\proclaim
\let \endpro=\endproclaim

\begin{document}
 \centerline {\bf
The residue-counts of $x^2+a/x$ modulo a prime}
\par\q\newline
\centerline{Zhi-Hong Sun}\newline
\centerline{School of Mathematics
and Statistics}
\centerline{Huaiyin Normal University}
\centerline{Huaian, Jiangsu 223300, P.R. China} \centerline{Email:
zhsun@hytc.edu.cn} \centerline{Homepage:
http://maths.hytc.edu.cn/szh1.htm}
 \abstract{\par For a prime $p>3$ and $a\in \Bbb Z$ with $p\nmid a$ let
$V_p(x^2+\f ax)$ be the residue-counts of $x^2+\f ax$ modulo $p$ as $x$ runs over $1,2,\ldots,p-1$. In this paper, we obtain an explicit formula for $V_p(x^2+\f ax)$,
which is concerned with cubic residues and binary quadratic forms.
 \par\q
 }
 \endabstract

\section*{1. Introduction}
\par\q\ 
For an odd prime $p$ let $\Bbb Z_p=\{0,1,\ldots,p-1\}$ and $\Bbb Z_p^*=\{1,2,\ldots,p-1\}$,
and let $R_p$ be the set of those rational numbers whose denominators are not divisible by $p$. For an odd prime $p$ and $a\in R_p$ let $\sls ap$ be the Legendre symbol, and let $V_p(f(x))$ be the residue-counts of $f(x)$ modulo $p$. Namely,
$V_p(f(x))$ is the number of $b\in\Bbb Z_p$ such that $f(x)\e b\mod p$ is solvable.
 \par Let $p>3$ be a prime. It is well known that $V_p(x^2)=\f{p+1}2$. In 1908, von Sterneck[6] showed that for $a_1,a_2,a_3\in\Bbb Z$ with
$a_1^2-3a_2\not\e 0\mod p$,
$$V_p(x^3+a_1x^2+a_2x+a_3)=\f{2p+\sls p3}3.$$
See also [2, Theorem 4.3]. Let $[x]$ be the greatest integer not exceeding $x$. Suppose  $b\in\Bbb Z$ with $p\nmid b$. For $p\e 2\mod 3$, in [4] the author proved that $V_p(x^4+bx)=[\f{5p+7}8]$. For $p=3k+1=A^2+3B^2$ with $A,B\in\Bbb Z$ and $3\mid A-1$, in [4] the author proved that
$$V_p\big(x^4+bx\big)=\cases \f 18\big(5p+9-6(-1)^{\f{p-1}{12}}\big)\q\t{if $12\mid p-1$ and $2b$ is a cubic residue of $p$,}
\\\f 18(5p+3\pm 6B)\q\t{if $12\mid p-1$ and
$(2b)^{\f{p-1}3}\e\f 12(-1\mp\f AB)\mod p$,}
\\\f 18\big(5p+7+6(-1)^{\f{p-7}{12}}-4A\big)
\q\t{if $12\mid p-7$ and $2b$ is a cubic residue of $p$,}
\\\f 18(5p+1+2A)\q\t{if $12\mid p-7$
and $2b$ is a cubic non-residue of $p$.}
\endcases$$
In a recent paper [5], Sun and Ye obtained a formula for $V_p(x^3+\f cx)$ for $c\in \Bbb Z$ with $p\nmid c$. In particular, for $p\e 3\mod 4$,
$$V_p\Big(x^3+\f cx\Big)=\cases
\f{5p+5}8+\f 18\big(\frac{-c}{p}\big)a_{24}(p) &\t{if $8\mid p-7$ and $\ls{3c}p=1$,}
\\\f{5p-3}8+\f 18\big(\frac{-c}{p}\big)a_{24}(p) &\t{otherwise,}
\endcases$$
where $a_{24}(n)$ is given by
$$q\prod_{k=1}^{\infty}(1-q^{2k})(1-q^{4k})(1-q^{6k})(1-q^{12k})
=\sum_{n=1}^{\infty}a_{24}(n)q^n\quad(|q|<1).$$
\par Motivated by the above work, in this paper we use the theory of cubic congruences in [2,3] to evaluate $V_p(x^2+\f ax)$, where $p>3$ is a prime and $a\in R_p$ with $a\not\e 0\mod p$.
In particular, we show that for $p=3k+1=A^2+3B^2$ with $A,B\in\Bbb Z$ and $3\mid A-1$,
$$V_p\Big(x^2+\f{2a}x\Big)=\cases
\f{2p-1+2A}3&\t{if $a^{\f{p-1}3}\e 1\mod p$,}
\\\f{2p-1-A\pm 3B}3&\t{if $a^{\f{p-1}3}\e \f {-1\mp A/B}2\mod p$}
\endcases$$
and so $a$ is a cubic residue of $p$ if and only if
 $V_p(x^2+\f {2a}x)=\f{2p-1+2A}3$. As consequences,
$$A=\f 12\Big(3V_p\Big(x^2+\f 2x\Big)+1\Big)-p=\f 12\Big(3V_p\Big(x+\f 4{x^2}\Big)+1\Big)-p$$
and
$$L=2p-1-3V_p\Big(x^2+\f 1x\Big)=2p-1-3V_p\Big(x+\f 1{x^2}\Big),$$
where $L$ is given by $4p=L^2+27M^2\ (L,M\in\Bbb Z)$ and $L\e 1\mod 3$. We remark that Jacobi
proved the remarkable congruences:
$$A\e \f 12\b{\f{p-1}2}{\f{p-1}6}\mod p\q\t{and}\q L\e -\b{\f{2(p-1)}3}{\f{p-1}3}\mod p.$$
\section*{2. Main results}
\par
The discriminant of
the cubic polynomial $x^3+a_1x^2+a_2x+a_3$ is given by $$D=a_1^2a_2^2-4a_2^3-4a_1^3a_3-27a_3^2+18a_1a_2a_3.\tag 2.1$$ Let $p>3$ be a prime, and let
$N_p(f(x))$
denote the number of
solutions of the congruence $f(x)\e 0\mod p$.

\par {\bf Lemma 2.1 ([3, Theorem 1.1, Lemmas 2.3 and 4.1])} If $p>3$ is a prime, $a_1,a_2,a_3\in \Bbb Z$ and $D$ is given by (2.1), then
$$N_p(x^3+a_1x^2+a_2x+a_3)=\cases 0\ \t{or}\ 3&\t{if $\sls Dp=1$},
\\3&\t{if $\sls Dp=0$,}\\ 1  &\t{if $\sls Dp=-1$}.\endcases$$
\par Lemma 2.1 is actually attributed to Stickelberger, Dickson and Skolem. See the references in [3].
\par Set $\omega=(-1+\sqrt{-3})/2$.  Let $p>3$ be a prime, $k\in R_p$ and $k^2+3\not\e 0\mod p$, and let $\big(\f{k+1+2\omega}p\big)_3$ be the cubic Jacobi symbol defined in [2]. Following [2] define
$$C_0(p)=\Big|\Big\{k\in R_p: k^2+3\not\e 0\mod p,\ \Big(\f{k+1+2\omega}p\Big)_3=1\Big\}\Big|.$$
By [2, Corollary 3.3],
$$k\in C_0(p)\iff k\e \f{x^3-9x}{3x^2-3}\mod p\ \t{for $x\in\{0,2,3,\ldots,p-2\}$ with $p\nmid x^2+3$}.\tag 2.2$$
Since $k\in C_0(p)$ implies $-k\in C_0(p)$ and $\big(\f{x^3-9x}{3x^2-3}\big)^2+3=\f{(x^2+3)^3}{9(x^2-1)^2}$, from (2.2) we deduce that
$$k\in C_0(p)\iff 9(k^2+3)\e \f{(x^2+3)^3}{(x^2-1)^2}\mod p\ \t{for $x\in H_p$},\tag 2.3$$
where
$$H_p=\Big\{x:\ x\in\big\{0,2,3,\ldots,\f{p-1}2\big\}\ \t{and}\ x^2\not\e -3\mod p\Big\}.$$

\pro{Lemma 2.2} Let $p>3$ be a prime, $t\in R_p$ and $t\not\e 0\mod p$. If $t\e \f{(x^2+3)^3}{(x^2-1)^2}\mod p$ is solvable for $x\in H_p$, then
$$\Big|\Big\{x\in H_p:\ \f{(x^2+3)^3}{(x^2-1)^2}\e t\mod p\Big\}\Big|
=\cases 3&\t{if $t\not\e 27\mod p$,}
\\ 2&\t{if $t\e 27\mod p$.}
\endcases$$
\endpro
{\it Proof}. Suppose $x_0\in H_p$, $t\e \f{(x_0^2+3)^3}{(x_0^2-1)^2}\mod p$ and
$c\in\Bbb Z_p$ with $c\e \f{x_0^3-9x_0}{3(x_0^2-1)}\mod p$. Then
$$t\e \f{(x_0^2+3)^3}{(x_0^2-1)^2}=9\Ls{x_0^3-9x_0}{3(x_0^2-1)}^2+27\e 9c^2+27\mod p.$$
Clearly,  the discriminant of $x^3-9x-3c(x^2-1)$ is $(18(c^2+3))^2$ and $c^2+3\not\e 0\mod p$. Since $x_0^3-9x_0-3c(x_0^2-1)\e 0\mod p$, from Lemma 2.1 we see that
    $x^3-9x-3c(x^2-1)\e 0\mod p$ has three distinct solutions and so $c\e \f{x^3-9x}{3x^2-3}\mod p$ has three distinct solutions. Note that $\f{(-x)^3-9(-x)}{3(-x)^2-3}=-\f{x^3-9x}{3x^2-3}$ and $1^2,2^2,\ldots\ls{p-1}2^2$ are pairwise distinct modulo $p$.
     For $t\not\e 27\mod p$ and so $c\not\e 0\mod p$, there are exactly three $x\in\{0,2,3,\ldots,\f{p-1}2\}$ such that $t\e 9(c^2+3)\e 9\ls{x^3-9x}{3x^2-3}^2+27=\f{(x^2+3)^3}{(x^2-1)^2}\mod p$.
     Recall that $t\not\e 0\mod p$. Thus, for $t\not\e 27\mod p$ there are exactly three
      $x\in H_p$ such that $t\e \f{(x^2+3)^3}{(x^2-1)^2}\mod p$.
   For $t\e 27\mod p$ and so $c\e 0\mod p$, clearly $c\e \f{x^3-9x}{3x^2-3}\mod p$ has three solutions $x\e 0,3,-3\mod p$. Since $3^2=(-3)^2$, there are exactly two $x\in\{0,2,3,\ldots,\f{p-1}2\}$ such that $9c^2+27\e 9\ls{x^3-9x}{3x^2-3}^2+27=\f{(x^2+3)^3}{(x^2-1)^2}\mod p$. Therefore, there are exactly two
      $x\in H_p$ such that $t\e \f{(x^2+3)^3}{(x^2-1)^2}\mod p$. The proof is now complete.
\vskip 0.2cm
\pro{Lemma 2.3} Let $p>3$ be a prime, $m\in R_p$ and $m\not\e 0\mod p$.
Then
$$\align&\Ls mp\sum_{y=1}^{p-1}\Ls{y^3+m}p
\\&=\cases -1\q\qq\qq\;\,\t{if $p\e 2\mod 3$,}
\\ -1-2A\q\qq\t{if $p=3k+1=A^2+3B^2(A,B\in\Bbb Z)$, $3\mid A-1$} \\\qq\qq\qq\qq\t{and $m^{\f{p-1}3}\e 1\mod p$,}
\\-1+A\mp 3B\q\t{if $p=3k+1=A^2+3B^2(A,B\in\Bbb Z)$, $3\mid A-1$}\\\qq\qq\qq\qq \t{and $m^{\f{p-1}3}\e \f {-1\pm A/B}2\mod p$.}
\endcases\endalign$$\endpro
\par Lemma 2.3 is a known result. See [1, pp.195-196] and [4,(2.9)]. The sum 
$\sum_{y=1}^{p-1}\ls{y^3+m}p$ is called cubic Jacobsthal sum. For $p\e 2\mod 3$
and $c\in\Bbb Z$, $y^3\e c\mod p$ has a unique solution. Thus, for $p\e 2\mod 3$,
$$\Ls mp\sum_{y=1}^{p-1}\Ls{y^3+m}p=\Ls mp\sum_{x=1}^{p-1}\Ls{x+m}p=\Ls mp\sum_{x=0}^{p-1}\Ls{x+m}p-1=-1.$$

\pro{Theorem 2.1} Let $p>3$ be a prime and $a\in R_p$ with $a\not\e 0\mod p$.
\par $(\rm i)$ If $p\e 2\mod 3$, then $V_p(x^2+\f ax)=\f{2p-1}3$.
\par $(\rm ii)$ If $p\e 1\mod 3$ and so $p=A^2+3B^2$ with $A,B\in\Bbb Z$ and $A\e 1\mod 3$, then
$$V_p\Big(x^2+\f ax\Big)=\cases
\f{2p-1+2A}3&\t{if $(2a^2)^{\f{p-1}3}\e 1\mod p$,}
\\\f{2p-1-A\pm 3B}3&\t{if $(2a^2)^{\f{p-1}3}\e \f {-1\pm A/B}2\mod p$.}
\endcases$$
\endpro
{\it Proof.}
Set
$$\delta_p(c)=\cases 1&\t{if $c$ is a cubic residue of $p$,}
\\0&\t{if $c$ is a cubic non-residue of $p$.}\endcases$$
For $m\in\Bbb Z$ and $p\e 2\mod 3$, the congruence $x^3\e m\mod p$ has a unique solution.
If $p\e 1\mod 3$ and $x^3\e m\mod p$ is solvable, then $x^3\e m\mod p$ has three solutions.
Since
\begin{align*} &\Big|\Big\{b\in\Bbb Z_p:\ \Ls{4b^3-27a^2}p=1\Big\}\Big|
+\Big|\Big\{b\in\Bbb Z_p:\ \Ls{4b^3-27a^2}p=-1\Big\}\Big|
\\&=p-\Big|\Big\{b\in\Bbb Z_p:\ 4b^3\e 27a^2\mod p\Big\}\Big|=p-\Big|\Big\{c\in\Bbb Z_p:\ c^3\e 2a^2\mod p\Big\}\Big|
\\&=p-\Big(2+\Ls p3\Big)\delta_p(2a^2)
\end{align*}
and
\begin{align*} &\Big|\Big\{b\in\Bbb Z_p:\ \Ls{4b^3-27a^2}p=1\Big\}\Big|
-\Big|\Big\{b\in\Bbb Z_p:\ \Ls{4b^3-27a^2}p=-1\Big\}\Big|
=\sum_{b\in\Bbb Z_p}\Ls{4b^3-27a^2}p,
\end{align*}
we see that
$$\Big|\Big\{b\in\Bbb Z_p:\ \Ls{4b^3-27a^2}p=-1\Big\}\Big|=\f 12
\Big(p-\sum_{b\in\Bbb Z_p}\Ls{4b^3-27a^2}p-\Big(2+\Ls p3\Big)\delta_p(2a^2)\Big).$$
It is clear that
\begin{align*}\sum_{b\in\Bbb Z_p}\Ls{4b^3-27a^2}p&=\sum_{y\in\Bbb Z_p}\Ls{4(-\f 32y)^3-27a^2}p=\Ls{-6}p
\sum_{y\in\Bbb Z_p}\Ls{y^3+2a^2}p\\&=\Ls{-3}p+\Ls{-6}p
\sum_{y=1}^{p-1}\Ls{y^3+2a^2}p.\end{align*}
Thus,
$$\aligned&\Big|\Big\{b\in\Bbb Z_p:\ \Ls{4b^3-27a^2}p=-1\Big\}\Big|
\\&=\f 12
\Big(p-\Ls{-3}p-\Ls{-6}p\sum_{y=1}^{p-1}\Ls{y^3+2a^2}p-\Big(2+\Ls p3\Big)\delta_p(2a^2)\Big).\endaligned\tag 2.4$$
 Observe that the discriminant of $x^3-bx+a$ is $4b^3-27a^2$. By Lemma 2.1, for $\sls{4b^3-27a^2}p=-1$ the congruence $x^3-bx+a\e 0\mod p$ has a unique solution. Thus, from (2.4) we have
    $$\aligned&\Big|\Big\{b\in\Bbb Z_p:\ \Ls{4b^3-27a^2}p=-1, \ x^2+\f ax\e b\mod p\ \t{is solvable}\Big\}\Big|
\\&=\f 12
\Big(p-\Ls p3-\Ls{-6}p\sum_{y=1}^{p-1}\Ls{y^3+2a^2}p-\Big(2+\Ls p3\Big)\delta_p(2a^2)\Big).\endaligned\tag 2.5$$
\par
  For $b,k\in\Bbb Z$ such that $b\not\e 0\mod p$, $\f{4b^3-27a^2}{9a^2}\e k^2\mod p$ and so $\f{4b^3}{a^2}\e 9(k^2+3)\mod p$, substituting $a,b$ with $\f b3,-\f{3a}b$ in [2, Theorem 4.1] we deduce that $x^3-bx+a\e 0\mod p$ is solvable if and only if $k\in C_0(p)$.
  Now, appealing to (2.3) we see that $x^3-bx+a\e 0\mod p$ is solvable if and only if
  $\f{4b^3}{a^2}\e \f{(x^2+3)^3}{(x^2-1)^2}\mod p$ for some $x\in H_p$. Hence,
      \begin{align*} &\Big|\Big\{b\in\Bbb Z_p^*:\ \Ls{4b^3-27a^2}p\in\{0,1\},\
  x^2+\f ax\e b\mod p\ \t{is solvable}\Big\}\Big|
  \\&=\Big|\Big\{b\in\Bbb Z_p^*:\ \Ls{4b^3-27a^2}p\in\{0,1\},\
  x^3-bx+a\e 0\mod p\ \t{is solvable}\Big\}\Big|
    \\&=\Big|\Big\{b\in\Bbb Z_p^*:\ \f{4b^3}{a^2}\e \f{(x^2+3)^3}{(x^2-1)^2}\mod p\ \t{for $x\in H_p$}\Big\}\Big|
        \\&=\Big|\Big\{b\in\Bbb Z_p^*:\ b^3\e 2a^2(x^2-1)\Ls{x^2+3}{2(x^2-1)}^3\mod p\ \t{for $x\in H_p$}\Big\}\Big|.
        \end{align*}
        If $2a^2(x^2-1)\e y^3\mod p$ for  $x\in H_p$ and $y\in\Bbb Z_p$, then for $b\in\Bbb Z_p^*$ with $b\e \f{x^2+3}{2(x^2-1)}y\mod p$ we have
    $b^3\e 2a^2(x^2-1)\ls{x^2+3}{2(x^2-1)}^3\mod p$. Conversely, if $b^3\e 2a^2(x^2-1)\ls{x^2+3}{2(x^2-1)}^3\mod p$ for $b\in\Bbb Z_p^*$ and $x\in H_p$, then $2a^2(x^2-1)\e y^3\mod p$ is solvable. Suppose $2a^2(x^2-1)\e y_i^3\mod p$ and $b_i\in\Bbb Z_p$ with $b_i\e\f{x^2+3}{2(x^2-1)}y_i\mod p$ for $i=1,\ldots,2+\ls p3$, where $y_i\in\Bbb Z_p$ and $y_1,\ldots,y_{2+\sls p3}$ are distinct. Then
    $b_i^3\e 2a^2(x^2-1)\ls{x^2+3}{2(x^2-1)}^3\e b^3\mod p$ for $i=1,\ldots,2+\sls p3$ and so $b\in\{b_1,\ldots,b_{2+\sls p3}\}$. Also,
    $$\big|\big\{b\in\Bbb Z_p^*:\ 4b^3\e 27a^2\mod p\big\}\big|=\Big(2+\Ls p3\Big)\delta_p(2a^2).$$
     \par Now, from the above and Lemma 2.2 we deduce that
    $$\align &3\Big|\Big\{b\in\Bbb Z_p^*:\ b^3\e 2a^2(x^2-1)\Ls{x^2+3}{2(x^2-1)}^3\mod p\ \t{for $x\in H_p$}\Big\}\Big|\\&=\Big|\Big\{(x,y)\in H_p\times \Bbb Z_p^*:\ y^3\e 2a^2(x^2-1)\mod p\Big\}\Big|+\Big(2+\Ls p3\Big)\delta_p(2a^2)
    \\&=\Big|\Big\{(x,y)\in\big\{0,1,\ldots,\f{p-1}2\big\}\times\Bbb Z_p:\ y^3\e 2a^2(x^2-1)\mod p\Big\}\Big|
    \\&\q-1-\f 32\Big(1+\Ls p3\Big)\delta_p(a^2)+\Big(2+\Ls p3\Big)\delta_p(2a^2)
    \\&=\f 12\Big|\Big\{(x,y)\in\Bbb Z_p\times\Bbb Z_p:\ y^3\e 2a^2(x^2-1)\mod p\Big\}\Big|
    \\&\q+\f 32\Big(2+\Ls p3\Big)\delta_p(2a^2)
    -1-\f 32\Big(1+\Ls p3\Big)\delta_p(a^2)
    \endalign$$
and so
\begin{align*} &\Big|\Big\{b\in\Bbb Z_p:\ \Ls{4b^3-27a^2}p\in\{0,1\},\
  x^2+\f ax\e b\mod p\ \t{is solvable}\Big\}\Big|
  \\&=\Big|\Big\{b\in\Bbb Z_p^*:\ \Ls{4b^3-27a^2}p\in\{0,1\},\
  x^2+\f ax\e b\mod p\ \t{is solvable}\Big\}\Big|+\f {1+\sls p3}2\delta_p(a)
  \\&=\f 16\Big|\Big\{(x,y)\in\Bbb Z_p\times\Bbb Z_p:\ y^3\e 2a^2(x^2-1)\mod p\Big\}\Big|
    +\f 12\Big(2+\Ls p3\Big)\delta_p(2a^2)-\f 13.
\end{align*}
On the other hand,
$$\aligned
&\big|\big\{(x,y)\in\Bbb Z_p\times \Bbb Z_p:\ 2a^2(x^2-1)\e y^3\mod p\big\}\big|
\\&=\Big|\Big\{(x,y)\in\Bbb Z_p\times \Bbb Z_p:\ x^2\e \f{y^3}{2a^2}+1\mod p\Big\}\Big|
\\&=\sum\Sb y=0\\ \sls{1+y^3/(2a^2)}p=1\endSb^{p-1}2+\sum\Sb y=0\\ \sls{1+y^3/(2a^2)}p=0\endSb^{p-1}1
=\sum\Sb y=0\\ \sls{1+y^3/(2a^2)}p=1\endSb^{p-1}1+p-\sum\Sb y=0\\ \sls{1+y^3/(2a^2)}p=-1\endSb^{p-1}1
\\&=p+\sum_{y=0}^{p-1}\Ls {1+y^3/(2a^2)}p
=p+1+\Ls 2p\sum_{y=1}^{p-1}\Ls {y^3+2a^2}p.
\endaligned$$
Thus,
\begin{align*} &\Big|\Big\{b\in\Bbb Z_p:\ \Ls{4b^3-27a^2}p\in\{0,1\},\
  x^2+\f ax\e b\mod p\ \t{is solvable}\Big\}\Big|
  \\&=\f 16\Big(p+1+\Ls 2p\sum_{y=1}^{p-1}\Ls {y^3+2a^2}p\Big)
  +\f 12\Big(2+\Ls p3\Big)\delta_p(2a^2) -\f 13.
    \end{align*}
    This together with (2.5) yields
   $$\align V_p\Big(x^2+\f ax\Big)&=\f{p-\sls p3}2-\f 12\Ls{-6}p\sum_{y=1}^{p-1}\Ls{y^3+2a^2}p-\f{2+\sls p3}2\delta_p(2a^2)
  \\&\q+\f {p+1}6+\f 16\Ls 2p\sum_{y=1}^{p-1}\Ls {y^3+2a^2}p
  +\f 12\Big(2+\Ls p3\Big)\delta_p(2a^2)
    -\f 13\endalign$$
    and so
    $$V_p\Big(x^2+\f ax\Big)=\f {2(p-\sls p3)}3+\f{\sls p3-1}6
    +\f 16\Big(1-3\Ls p3\Big)\Ls 2p\sum_{y=1}^{p-1}\Ls{y^3+2a^2}p.\tag 2.6$$
Now, applying Lemma 2.3 we derive that
$$V_p\Big(x^2+\f ax\Big)=\cases \f{2(p-1)}3-\f 13(-1-2A)=\f{2p-1+2A}3
\\\qq\qq\qq\t{if $(2a^2)^{\f{p-1}3}\e 1\mod p$,}
\\\f{2(p-1)}3-\f 13(-1+A\mp 3B)=\f{2p-1-A\pm 3B}3
\\\qq\qq\qq\t{if $(2a^2)^{\f{p-1}3}\e \f {-1\pm A/B}2\mod p$.}
\endcases$$
This completes the proof.
\vskip 0.2cm
\pro{Corollary 2.1} Let $p$ be a prime of the form $3k+1$
and so $p=A^2+3B^2$ with $A,B\in\Bbb Z$ and $A\e 1\mod 3$,
and let $a\in R_p$ with $a\not\e 0\mod p$.
Then
$$V_p\Big(x^2+\f{2a}x\Big)=\cases
\f{2p-1+2A}3&\t{if $a^{\f{p-1}3}\e 1\mod p$,}
\\\f{2p-1-A\pm 3B}3&\t{if $a^{\f{p-1}3}\e \f {-1\mp A/B}2\mod p$}
\endcases$$
and so
$a$ is a cubic residue of $p$ if and only if
 $V_p(x^2+\f {2a}x)=\f{2p-1+2A}3$.
\endpro
{\it Proof}. Replacing $a$ with $2a$ in Theorem 2.1 and noting that
$a^{-\f{p-1}3}\e \f {-1\mp A/B}2\mod p$ implies $a^{\f{p-1}3}\e \f {-1\pm A/B}2\mod p$
we deduce the formula for $V_p(x^2+\f {2a}x)$.
Since $p$ is a prime we see that $A^2\not=B^2$ and so $A\not=\pm B$. Hence, if $a$ is a cubic non-residue modulo $p$, from the above we have
 $V_p(x^2+\f {2a}x)=\f{2p-1+2A-3A\pm 3B}3\not=\f{2p-1+2A}3$.
 Thus the result follows.
\vskip 0.2cm
\pro{Corollary 2.2} Let $p>3$ be a prime and $a\in R_p$ with $a\not\e 0\mod p$.
\par $(\rm i)$ If $p\e 2\mod 3$, then $V_p(x+\f a{2x^2})=\f{2p-1}3$.
\par $(\rm ii)$ If $p\e 1\mod 3$ and so $p=A^2+3B^2$ with $A,B\in\Bbb Z$ and $A\e 1\mod 3$, then
$$V_p\Big(x+\f a{2x^2}\Big)=\cases
\f{2p-1+2A}3&\t{if $a^{\f{p-1}3}\e 1\mod p$,}
\\\f{2p-1-A\pm 3B}3&\t{if $a^{\f{p-1}3}\e \f {-1\pm  A/B}2\mod p$.}
\endcases$$
\endpro
{\it Proof.} For $x\in\Bbb Z_p^*$, it is well known that there exists a unique $x'\in\Bbb Z_p^*$ such that $xx'\e 1\mod p$ and so $\f 1{{x'}^2}+cx'\e x^2+\f cx\mod p$ for $c\in R_p$. Thus
$$V_p\Big(x+\f a{2x^2}\Big)=V_p\Big(\f 1x+\f {ax^2}2\Big)=V_p\Big(x^2+\f 2{ax}\Big).\tag 2.7$$
Now, replacing $a$ with $\f 1a$ in Corollary 2.1 yields the result.
\vskip 0.2cm

\pro{Corollary 2.3} Let $p$ be a prime of the form $3k+1$ and so $4p=L^2+27M^2$ with $L,M\in\Bbb Z$ and $L\e 1\mod 3$. Then
$$L=2p-1-3V_p\Big(x^2+\f 1x\Big)=2p-1-3V_p\Big(x+\f 1{x^2}\Big).$$
\endpro
{\it Proof.}  Suppose $p=A^2+3B^2$ with $A,B\in\Bbb Z$ and $3\mid A-1$. From [4, (2.10)-(2.12)],
$$L=\cases -2A &\t{if $2^{\f{p-1}3}\e 1\mod p$,}
\\A+3B&\t{if $3\mid B-1$ and so $2^{\f{p-1}3}\e \f 12(-1-\f AB)\mod p$.}
\endcases$$
Now, taking $a=\f 12$ in Corollary 2.1 yields $V_p(x^2+\f 1x)=\f{2p-1-L}3$. To complete the proof, we note that $V_p(x^2+\f 1x)=V_p(x+\f 1{x^2})$ by (2.7).
\vskip 0.2cm

\pro{Corollary 2.4} Let $p$ be a prime of the form $3k+1$ and so $p=A^2+3B^2$ with $A,B\in\Bbb Z$ and $A\e 1\mod 3$. Then
$$V_p\Big(x^2+\f 4x\Big)=V_p\Big(x+\f 2{x^2}\Big)=\cases \f{2p-1+2A}3&\t{if $3\mid B$,}
\\\f{2p-1-A+3\sls B3B}3&\t{if $3\nmid B$.}
\endcases$$
\endpro
{\it Proof.} It is well known that $2^{\f{p-1}3}\e 1\mod p$ if and only if $3\mid B$. When  $2^{\f{p-1}3}\not\e 1\mod p$, we have $2^{\f{p-1}3}\e \f 12(-1-\sls B3\f AB)\mod p$ by [4, (2.12)]. Now, taking $a=2$ in Corollary 2.1 and
$a=4$ in Corollary 2.2 yields the result.
\vskip 0.2cm


\end{document}